\documentclass{article}

\begin{document}

\title{Nahm's equations  and free boundary problems }
\author{S. K. Donaldson\\Imperial College, London}
\maketitle
\ \ \ \ \ \ \ \  {\it Dedicated to Nigel Hitchin, with gratitude and affection}

\newcommand{\bR}{{\bf R}}
\newcommand{\bZ}{{\bf Z}}
\newcommand{\bC}{{\bf C}}
\newcommand{\curl}{{\rm curl}\ }
\newcommand{\diver}{{\rm div}\ }
\newcommand{\grad}{{\rm grad}\ }
\newcommand{\ux}{\underline{x}}
\newcommand{\up}{\underline{p}}
\newcommand{\uq}{\underline{q}}
\newtheorem{thm}{Theorem}
\newtheorem{lem}{Lemma}
\newtheorem{cor}{Corollary}
\newtheorem{prop}{Proposition}
\newtheorem{Question}{Question}

\section{Introduction}
In \cite{kn:D1}, following up work of Hitchin \cite{kn:H}, the author found
it useful to express {\it Nahm's equations}, for a matrix group,  in terms of the motion of a particle
in a Riemannian symmetric space, subject to  a potential field. This point of view
lead readily to an elementary existence theorem for solutions of Nahm's equation, corresponding to particle
paths with prescribed end points. The original motivation for this article
is the question of formulating an analogous theory for the  Nahm equations associated
to the infinite-dimensional Lie group of area-preserving diffeomorphisms
of a surface--in the spirit of \cite{kn:D2}. We will see that this can be done, and that a form of the appropriate
existence theorem holds---essentially a special case of a result of
Chen. However the main focus of the article is not on existence proofs but
on the various formulations of the problem, and connections between them.
In these developments, one finds that the natural context 
is rather more general than the original question, so we will start out of
a different tack, and return to Nahm's equations in Section 5.

Consider the following set-up in Euclidean space  $\bR^{3}$, in which
we take coordinates $(x_{1},x_{2}, z)$---thinking  of $z$ as the
vertical direction. (We will use the notation $\frac{\partial}{\partial x_{i}}=
\partial_{i}, \frac{\partial}{\partial z}= \partial_{z}$.) Suppose we have a strictly
positive function
$H(x_{1}, x_{2})$. This defines a domain 
$$\Omega_{H}= \{ (x_{1},x_{2}, z)\ : \ 0< z< H(x_{1},x_{2})\},
$$
whose boundary has two components $\{ z=0\}$ and $ \{ z=H\}$. We consider the Dirichlet
problem for the standard Laplacian: to find a harmonic function
$\theta$ on $\Omega_{H}$ with $\theta=0$ on $\{z=0\}$ and $\theta=1$ on $\{ z=H\}$.
To set up this problem precisely, let us assume that the data $H$ is
 $\bZ^{2}$-periodic on $\bR^{2}$, and seek a $\bZ^{2}$-periodic solution $\theta$. 
 Now we have a unique solution
$\theta$ to our Dirichlet problem. Consider the flux of the gradient of $\theta$ through the boundary $\{ z=H\}$. This defines another
function $\rho$ on $\bR^{2}$. To be precise, if  $\iota_{H}$ is the obvious
map from $\bR^{2}$ to the boundary $\{ z=H\}$ then the flux is defined by
\begin{equation}  \iota_{H}^{*}(* d\theta)= \rho dx_{1} dx_{2}. \label{eq:flux1} \end{equation}
 Explicitly
$$  \rho= \partial_{z} \theta - ( \partial_{1} \theta \partial_{1} H+ \partial_{2} \theta \partial_{2}
H), $$
with the right hand side evaluated at $ (x_{1}, x_{2}, H(x_{2}, x_{2}))$.
By the maximum principle, $\rho$ is a positive
function, since the normal derivative of $\theta$ is positive in the positive
$z$ direction on $\{ z=H \}$. We consider the following free boundary problem: given a positive periodic
function $\rho$ does it arise from some periodic $H$, and is $H$ unique?

One can gain some physical intuition for this question by supposing that that the lower half-space
$\{z\leq 0\}$ represents a body with an infinite specific heat capacity fixed
at temperature $0$ and  $\Omega_{H}$
corresponds to a layer of ice covering this body. We choose units so that
the melting temperature
of the ice is $1$.  Sunlight  shines vertically downwards
onto  the upper surface $\{ z=H\}$ of the ice, but with a variable
intensity so that heat is transmitted to the surface according to the density
function $\rho$. We suppose that the surface of the ice is  sprinkled
by rain, which will instantly freeze if the surface temperature of the ice is less
than $1$. We also suppose that a  wind blows across
the surface, instantly removing
any surface water. Then  we see that the solution to our free boundary problem represents
a static physical state, in which the upper surface of the ice is just at
freezing point, the lower surface is at the imposed sub-freezing temperature  and the heat generated by the given sunlight flows through the ice
without changing the temperature. Physical intuition suggests that there
should indeed be a unique solution.

We can express the free-boundary problem considered above as a special case
of another question. Suppose now that we have a pair of periodic functions
$H_{0}, H_{1}$ on $\bR^{2}$ with $H_{0}< H_{1}$. Then we have a domain
$ \Omega_{H_{0}, H_{1}} = \{  H_{0}(x_{1}, x_{2})< z
< H_{1}(x_{1}, x_{2})\}, $
with two boundary components.
Let $\theta$ be the  harmonic function in this domain equal to $0, 1$ on $\{ z=
H_{0}\}, \{ z=H_{1}\}$ respectively. Then we obtain  a pair of flux-functions $\rho_{0}, \rho_{1}$ as
before. By Gauss' Theorem, these  satisfy a constraint
\begin{equation}  \int_{[0,1]^{2}} \rho_{0}\ d\ux= \int_{[0,1]^{2}} \rho_{1} \ d\ux, \label{eq:gauss}\end{equation}
since $ [0,1]^{2}$ is a fundamental domain for the $\bZ^{2}$-action.
Obviously, if we replace $H_{0}, H_{1}$ by $H_{0}+ c, H_{1} + c$ for any
constant $c$ we get the same pair $\rho_{0}, \rho_{1}$. We ask: given $\rho_{0},
\rho_{1}$ satisfying the integral constraint (2), is there a corresponding
pair $(H_{0}, H_{1})$, and if so is the solution unique up to the
addition of a constant? A positive answer to this question implies a positive
answer to the previous one, by a simple reflection argument. (Given
$\rho$, as in the first problem, take $\rho_{0}=\rho_{1}=\rho/2$. Then uniqueness implies  that the solution has reflection symmetry about $\theta=1/2$ and
we get a solution to the first problem by changing $\theta$ to $2\theta-1$.)

Of course we can also imagine a physical problem corresponding to this second
question: for example a layer of ice in the region $\Omega_{H_{0}, H_{1}}$.
We can now vary the problem by supposing that in place of ice we have a horizontally
stratified material in which heat can only flow in the horizontal directions.
Thus the steady-state condition, for a temperature distribution $\theta(x_{1},
x_{2}, z)$ is
\begin{equation}  (\partial_{1}^{2} +\partial_{2}^{2})\theta =0. \label{eq:partialDelta} \end{equation}
We define flux-functions $\rho_{0}, \rho_{1}$ by pulling back the $2$-form
$$   \partial_{1} \theta\ dx_{2} dz - \partial_{2} \theta \ dx_{1} dz, $$
 and the integral constraint
(2) still holds. So we ask: given $\rho_{0}, \rho_{1}$ satisfying (2), is
there a pair $H_{0}, H_{1}$ and a function $\theta$ on $\Omega_{H_{0}, H_{1}}$,
equal to $0,1$ on the two boundary components, which has these fluxes, and
is the solution essentially unique? (In this case
one has to relax the condition on the domain to $H_{0}\leq H_{1}$.)

It is natural to extend these questions to a general compact oriented Riemannian manifold
$X$ (which would be the flat torus $\bR^{2}/\bZ^{2}$ in the discussion above).
Write $d\mu$ for the Riemannian volume form on $X$.
We fix a real parameter $\epsilon\geq 0$ and define a  map $*_{\epsilon}$
from $T^{*}(X\times \bR)$ to $\Lambda^{n}T^{*}(X\times \bR)$ by 
$$  *_{\epsilon} dz= \epsilon d\mu\ \ , \ \ *_{\epsilon} \alpha= (*_{X}\alpha) dz,
$$
for $\alpha \in T^{*}X$, where $*_{X}$ is the usual Hodge $*$-operator on
$X$. Then for a function $\theta$ on a domain in $X\times \bR$ 
$$  d*_{\epsilon} d\theta = (\Delta_{\epsilon} \theta) dz d\mu, $$
where $$ \Delta_{\epsilon} \theta = (- \epsilon \partial_{z}^{2}+ \Delta_{X}) \theta,
$$
with $\Delta_{X}$ the standard Laplace operator on $X$. (We use the sign
convention that $\Delta_{X}$ is a positive operator, so when $\epsilon=1$
our $\Delta_{\epsilon}$ is the standard Laplace operator on $X\times \bR$.) If $\theta$ is defined on a domain $\Omega_{H_{0}, H_{1}}$, as above, we
define the  flux $\rho_{i}$ on the boundary $\{ z=H_{i}\}$ by pulling back
$*_{\epsilon} d\theta$ just as before.
We consider a pair of functions $\rho_{0},
\rho_{1}>0$ with
$$ \int_{X} \rho_{0}\ d\mu= \int_{X} \rho_{1} \ d\mu=\int_{X} \ d\mu $$
and we ask
\begin{Question}
Is there a pair $H_{0}\leq H_{1}$ and a function $\theta$ on the set 
$\Omega_{H_{0}, H_{1}}\subset X\times \bR$ with $\theta=0,1$ on the hypersurfaces
$\{ z=H_{0}\}, \{ z=H_{1}\}$, with fluxes $\rho_{i}$ and with
$  \Delta_{\epsilon} \theta = 0 $? If so, is the solution essentially unique?
\end{Question}

For any $\epsilon>0$ the equation $\Delta_{\epsilon} \theta=0$ can be transformed
into the standard Laplace equation on the product, by rescaling the $z$ variable.
When $\epsilon=0$ the equation has a very different character: it is not
elliptic and we obviously do not have automatic interior regularity with
respect to
$z$.

\section{An infinite-dimensional Riemannian manifold}

We now start in a different direction. Given our compact Riemannian manifold
$X$ we let ${\cal H}$ be the set of functions $\phi$ on $X$ such that
$1-\Delta_{X} \phi>0$. We make ${\cal H}$ into a Riemannian manifold,  defining
the norm of a tangent vector $\delta \phi$ at a point $\phi$ by
$$  \Vert \delta \phi \Vert^{2}_{\phi}= \int_{X} (\delta \phi)^{2}\ (1-\Delta_{X}
\phi) d\mu. $$
Thus  a path $\phi(t)$  in ${\cal H}$, parametrised by $t\in [0,1]$ say,
is simply a function on $X\times [0,1]$ and the \lq\lq energy'' of the path
is
\begin{equation} \frac{1}{2}  \int_{0}^{1} \int_{X}  \left( \frac{\partial\phi}{\partial t}\right)^{2} ( 1-\Delta_{X} \phi)\  d\mu\ dt. \label{eq:action}\end{equation} 
When $X$ is $2$-dimensional and orientable, this definition coincides with
the metric on the space of \lq\lq Kahler potentials'' discussed by Mabuchi
\cite{kn:M}, Semmes \cite{kn:S} and the author \cite{kn:D3}. The general
context in those references is a compact Kahler manifold: here we are considering
a different extension of the $2$-dimensional case, and we will see that some new features emerge. The account below follows the approach in \cite{kn:D3}
closely. 

It is  straightforward  to find the Euler-Lagrange equations associated
to the energy (3). These are
$$      \ddot{\phi} = \frac{\vert \nabla_{X} \dot{\phi}\vert^{2}}{1-\Delta_{X}
\phi}. $$
These equations define the geodesics in ${\cal H}$. We can read off the Levi-Civita
connection of the metric from this geodesic equation, as follows. Let
$\phi(t)$ be any path in ${\cal H}$ and $\psi(t)$ be another function on
$X\times [0,1]$, which we regard as a vector field along the path $\phi(t)$.
Then the covariant derivative of $\phi$ along the path is given by
\begin{equation} D_{t} \phi = \frac{d\phi}{dt} + (W_{t}, \nabla_{X}\psi),
\label{eq:covderiv}\end{equation}
where $$W_{t}= \frac{-1}{1-\Delta_{X} \phi} \nabla_{X} \dot{\phi}$$ and $(\ ,\ )$ is the Riemannian inner
product on tangent vectors to $X$. (We write $\nabla_{X}$, or sometimes just
$\nabla$, for the gradient operator on $X$, so $W_{t}$ is a vector field
on $X$.) This has an important consequence for
the {\it holonomy group} of the manifold ${\cal H}$. Observe that the tangent
space to ${\cal H}$ at a point $\phi$ is the space of functions on $X$ endowed
with the standard $L^{2}$ inner product associated to the measure
$$  d\mu_{\phi} = (1-\Delta \phi) d\mu_{0}. $$
So, in a general way, the parallel transport along a path from $\phi_{0}$
to $\phi_{1}$ should be an isometry from $L^{2}(X,d\mu_{\phi_{0}})$ to $L^{2}(X,d\mu_{\phi_{1}})$.
(Here we are ignoring the distinction between, for example, smooth functions
and $L^{2}$ functions.) What we see from equation (5) is that this isometry
is induced by a diffeomorphism $f:X\rightarrow X$ with
 \begin{equation} f^{*}(d\mu_{\phi_{1}})=d\mu_{\phi_{0}}. \label{eq:volid}\end{equation}
The diffeomorphism is obtained by integrating the time-dependendent vector
field $W_{t}$ and equation (6) follows from the identity
$$   {\cal L}_{W_{t}} d\mu_{\phi}= d* ( \frac{1}{1-\Delta_{X}} d\dot{\phi}*_{X}d\mu_{\phi})
= \Delta \dot{\phi}=-\frac{d}{dt} \mu_{\phi}. $$
(Here ${\cal L}$ denotes the Lie derivative on $X$.) We conclude that the
holonomy group of ${\cal H}$ is contained in the group ${\cal G}$ of volume-preserving
diffeomorphisms of $(X, d\mu_{0})$, regarded as a subgroup of the orthogonal
group of $L^{2}(X, d\mu_{0})$. (This can also be expressed by saying that
there is an obvious principal ${\cal G}$-bundle over ${\cal H}$ with the tangent
bundle as an associated vector bundle, and the Levi-Civita connection is
induced by a connection on this ${\cal G}$-bundle.)
 
  We now move on to discuss the curvature tensor of ${\cal H}$. Let $\phi$
  be a point of ${\cal H}$ and let $\alpha, \beta$ be tangent vectors to
  ${\cal H}$ at $\phi$---so $\alpha$ and $\beta$ are just functions on $X$.
  The curvature $R_{\alpha, \beta}$ should be a linear map from tangent vectors
  to tangent vectors: that is from functions on $X$ to functions on $X$.
  The discussion of the holonomy above tells us that this map must have the
  form
  \begin{equation}  R_{\alpha, \beta} (\psi) = (\nu_{\alpha, \beta}, \nabla \psi), \label{eq:curvid} \end{equation}
  for some vector field $\nu_{\alpha,\beta}$ on $X$, determined by
  $\phi, \alpha, \beta$. Moreover we know that we must have
  $$  {\cal L}_{\nu_{\alpha, \beta}} (d\mu_{\phi}) = 0. $$
  To identify this vector field we introduce some notation. For vector fields
  $v,w$ on $X$ we write $v\times w$ for the exterior product: a section of
  the bundle $\Lambda^{2} TX$. We define a differential operator
  $$  \curl: \Gamma(\Lambda^{2} TX) \rightarrow \Gamma(TX), $$
  to be the composite of the standard identification:
  $$  \Lambda^{2} TX \cong \Lambda^{n-2} T^{*} X, $$
    (using the Riemannian volume form $d\mu$), the exterior derivative
    $$ d: \Gamma(\Lambda^{n-2} T^{*} X\rightarrow \Lambda^{n-1} T^{*} X,
    $$
    and the standard identification
    $$   \Lambda^{n-1} T^{*} X\cong TX, $$
    (using the volume form $d\mu$ again).  Then we have
    \begin{thm}
    The curvature of ${\cal H}$ is given by (7) and the vector field
    $$  \nu_{\alpha, \beta} = \frac{1}{1-\Delta \phi} \curl(\frac{1}{1-\Delta
    \phi} \nabla \alpha \times  \nabla \beta). $$
    \end{thm}  

    \begin{cor}
    The manifold ${\cal H}$ has non-positive sectional curvature. 
    \end{cor}
    The sectional curvature corresponding to a pair of tangent vectors
    $\alpha, \beta$ at a point $\phi$ is
    $$  K_{\alpha, \beta}= \langle R_{\alpha, \beta} (\alpha), \beta\rangle.
    $$
    In our case this is
    $$  K_{\alpha, \beta}= \int_{X} (\nu_{\alpha, \beta}, \nabla \alpha) \beta (1-\Delta_{X}
    \phi) d\mu.$$
    Unwinding the algebraic identifications we used above, the integrand can be written in terms of differential forms
    as
    $$          \frac{1}{1-\Delta_{X}\phi} d\alpha \wedge d \left(\frac{1}{1-\Delta_{X}\phi}
    *(d\alpha \wedge d\beta)\right) \ \beta (1-\Delta_{X} \phi). $$
    So
    $$  K_{\alpha, \beta} =\int_{X} d\alpha \wedge d\left(\frac{1}{1-\Delta_{X}\phi} *(d\alpha \wedge
    d\beta)\right) \beta. $$
    Applying Stokes' Theorem this is
    $$ K_{\alpha, \beta}= -\int_{X}\frac{1}{1-\Delta_{X} \phi} d\alpha \wedge d\beta \wedge*(d\alpha \wedge d\beta)=- \int_{X}\frac{1}{1-\Delta_{X}\phi} \vert d\alpha
\wedge d\beta\vert^{2} \  d\mu\leq 0. $$

\

In the proof of Theorem 1 we will make use of two identities. For any pair
of vector fields $v,w$ and function $f$
\begin{equation} \curl (v\times w) = [v,w] + (\diver v) w - (\diver w) v \label{eq:id1}
\end{equation}
\begin{equation} \curl (f (v\times w)) = f \curl(v\times w)+ (v, \nabla f)
w-(w,\nabla f) v.\label{eq:id2} \end{equation}
We leave the verification as an exercise. (Considering geodesic coordinates
we see that it suffices to treat the case of Euclidean space. Our notation has be chosen to agree with standard notation in the case of
vector fields in $\bR^{3}$.)

To calculate the curvature we consider a $2$-parameter family $\phi(s,t)$
in ${\cal H}$, with a corresponding vector field $\psi(s,t)$  along the family.
Then we will compute the commutator
$  (D_{s} D_{t}- D_{t}D_{s}) \psi(s,t) $. Evaluating at $\phi=\phi(0,0)$ this
is $R_{\alpha, \beta}(\psi)$ where $\psi=\psi(0,0), \alpha=\partial_{s}\phi,
\beta=\partial_{t}\phi $.

Now we write $$D_{s}= \frac{\partial\ }{\partial s} + W_{s}\ , \ D_{t}= \frac{\partial
\ }{\partial t} + W_{t}, $$
where the vector fields $W_{s}, W_{t}$ are regarded as operators on the functions
on $X$. So $D_{s}D_{t}-D_{t}D_{s}$ is the operator given by the vector field
$$ \nu=    \frac{\partial W_{s}}{\partial t}- \frac{\partial W_{t}}{\partial s}
- [W_{s}, W_{t}], $$
and $\nu$ is exactly the vector field $\nu_{\alpha, \beta}$ we need to identify.
Recall that
$$   W_{s}= \frac{-\nabla \partial_{s}\phi}{1-\Delta \phi}\ , \  W_{t}= \frac{-\nabla \partial_{t}\phi}{1-\Delta
\phi}. $$
So
$$ \frac{\partial W_{s}}{\partial t}= \frac{-1}{1-\Delta \phi} \nabla\left(\frac{\partial^{2}
\phi}{\partial s \partial t}\right) + \frac{1}{(1-\Delta\phi)^{2}} \Delta \partial_{s}\phi
\nabla \partial_{t}\phi. $$
Evaluating at $s=t=0$ where $\partial_{s}\phi=\alpha,\partial_{t} \phi=\beta$ we have
$$  \frac{\partial W_{s}}{\partial t}- \frac{\partial W_{t}}{\partial s} =
\frac{1}{(1-\Delta \phi)^{2}} \left( \Delta \alpha \nabla \beta - \Delta
\beta \nabla \alpha\right). $$
Write $g$ for the function $(1-\Delta \phi)^{-1}$.
Combining with the Lie bracket term we obtain
$$  \nu_{\alpha, \beta} = [g \nabla \alpha, g \nabla \beta] + g^{2}( \Delta
\alpha \nabla \beta - \Delta \beta \nabla \alpha). $$

Now applying (8) we have
$$ [ g \nabla \alpha, g \nabla \beta] = \curl( g^{2} \nabla \alpha \times \nabla
\beta) + \diver(g \nabla \alpha) \nabla \beta- \diver( g \nabla \beta) \nabla
\alpha. $$
Applying (9) we have
$$   \curl( g^{2} \nabla \alpha \times \nabla \beta)= g\ \curl( g \nabla \alpha
\times \nabla \beta) + g( (\nabla g, \nabla \alpha)\nabla \beta - (\nabla
g, \nabla \beta) \nabla \alpha). $$
Since $$ \diver( g\nabla \alpha)= g \Delta \alpha - (\nabla g, \nabla \alpha)\
,\ 
\diver(g \nabla \beta)= g \Delta \beta-(\nabla g, \nabla \beta) $$
we see that $$ \nu_{\alpha, \beta}= g\curl( g \nabla \alpha \times \nabla \beta),
$$ as required.

In the case when $X$ has dimension $2$---as discussed in \cite{kn:D3}, \cite{kn:M},
\cite{kn:S}--- the space ${\cal H}$ is formally a symmetric space. This is
not true in general, since the curvature tensor is not preserved by the action
of the group ${\cal G}$.

We define a functional on ${\cal H}$ by

$$  V(\phi)= \int_{X} \phi\  d\mu. $$
This function is {\it convex} along geodesics in ${\cal H}$, since the geodesic
equation implies $\ddot{\phi}\geq
0$. Now introduce a real parameter $\epsilon\geq 0$ and consider the functional
on paths in ${\cal H}$:
\begin{equation}  E= \int \frac{1}{2} \vert \dot{\phi} \vert_{\phi}^{2} + \epsilon V(\phi) \ dt,
\label{eq:lagrangian} \end{equation}
corresponding to the motion of a particle in the potential $-\epsilon V$. The Euler-Lagrange equations are
\begin{equation}
\ddot{\phi}= \frac{\vert \nabla_{X} \dot{\phi}\vert^{2} + \epsilon}{1-\Delta_{X}
\phi}. \label{eq:phieq} \end{equation}

\section{Three equivalent problems}

In this section we will explain that there are three equivalent formulations
of the same PDE problem associated to a compact Riemannian manifold $X$. We have essentially encountered
two of these already.

\begin{itemize}

\item  The \lq\lq $\theta$ equation''.

\

This is the problem we set up in Section 1. We are given positive functions $\rho_{0},
\rho_{1}$ on $X$, with
\begin{equation} \int_{X} \rho_{i} \ d\mu = \int_{X} d\mu. \label{eq:normalint}
\end{equation}
We seek a domain $\Omega_{H_{0},H_{1}}\subset X\times
\bR$ defined by $H_{0}, H_{1}: X\rightarrow \bR$ and a function $\theta$
on $\Omega_{H_{0}, H_{1}}$,  equal to $0,1$ on the two boundary components,
with fluxes $\rho_{0}, \rho_{1}$ and satisfying the equation
$$  \Delta_{\epsilon}
\theta=0. $$

\

\item The \lq\lq $\Phi$ equation''

\

Here we are given $\phi_{0}, \phi_{1}$ on $X$, with $1-\Delta \phi_{i}>0$.
We seek a function $\Phi$ on $X\times [0,1]$, equal to $\phi_{0}, \phi_{1}$
on the two boundary components, with $1-\Delta \Phi>0$ for all $t$ and satisfying
the nonlinear equation
\begin{equation}   \frac{\partial^{2} \Phi}{\partial t^{2}} (1-\Delta_{X} \Phi)- \vert
\nabla \left( \frac{\partial \Phi}{\partial t}\right) \vert^{2} = \epsilon.
\label{eq:Phiequ} \end{equation}
As we have explained in Section 2, this is the same as finding a path in
the space ${\cal H}$, with end points $\phi_{0}, \phi_{1}$, corresponding to
the motion of a particle in the potential $-\epsilon V$. 

\

Now we introduce the  third problem.

\

\item  The \lq\lq $U$ equation''

We are given positive functions $\phi_{0}, \phi_{1}$, with $1-\Delta \phi_{i}>0,$ as 
above.  Define a function $L$ on $X\times \bR$ by
$$  L(x,z)= \max(\phi_{0}(x)-\phi_{1}(x) + z, 0). $$
We seek a $C^{1}$ function $U(x,z)$ on $X\times \bR$ with $U\geq L$ everywhere and satisfying the equation
   \begin{equation}   \Delta_{\epsilon} U= (1-\Delta \phi_{0})
    \label{eq:Ueq} \end{equation} on the open set $\Omega$ where $U>L$.
  \end{itemize}

   The equivalence of these three problems (assuming suitable regularity
   for the solutions in each case) arises from elementary, but not completely
   obvious, transformations. We describe these now.

   \begin{itemize}
   
  \item  $\theta$-equation \ $\Longrightarrow$\ $\phi$-equation
   
  \
   
  Suppose we have a solution $\theta$ on a domain $\Omega_{H_{0},H_{1}}$. Then $\partial_{z}\theta= \frac{\partial \theta}{\partial z}$
is positive on the boundary components of $\Omega_{H_{0}, H_{1}}$. The function
$\partial_{z}\theta$ satisfies the equation $ \Delta_{\epsilon} (\partial_{z}\theta) =0 $ and it follows from this that $\partial_{z}\theta$ is  positive throughout the domain. 
This implies that, for any $t\in [0,1]$, the set $\theta^{-1}(t)$ is the
graph of a smooth function 
 $h_{t}$ on $X$. By definition $h_{0}=H_{0}$ and $h_{1}=H_{1}$.
 We also write this function as $h(t,x)$
where convenient. 
For each fixed $t$ we can define a function $\rho_{t}$ on $X$ by the flux of
$ *_{\epsilon}d \theta$, just as before.

We claim that
\begin{equation} \frac{\partial \rho_{t}}{\partial t} = \Delta_{X} h_{t}\end{equation}

We show this by direct calculation (there are more conceptual, geometric arguments). For simplicity we treat the case when the metric on $X$ is locally
Euclidean, so $\Delta_{X}= - \sum \partial_{i}^{2}$ where $\partial_{i}=
\frac{\partial}{\partial x_{i}}$, for local coordinates $x_{i}$.  The identity
$$   \theta( x, h_{t} (x)) = t $$
implies that 
\begin{equation}  \partial_{i} \theta +\partial_{z} \theta \ \partial_{i} h= 0 \label{eq:id1}\end{equation}
and
\begin{equation} \partial_{z} \theta\ \partial_{t} h= 1. \label{eq:id2}\end{equation}
Now
$$ \Delta_{X} h_{t} = -\sum_{i} (\partial_{i} + (\partial_{i} h) \partial_{z}
)\partial_{i} h, $$ and this is
$$ \Delta_{X} h_{t}= -\sum_{i}(\partial_{i} - \frac{\partial_{i} \theta}{\partial_{z}
\theta} \partial_{z}) \left(- \frac{\partial_{i} \theta}{\partial_{z} \theta}\right) $$
which is
$$  -\sum_{i} \left( \frac{\partial_{i}^{2} \theta}{\partial_{z} \theta} - 2\frac{\partial_{i}
\theta \partial_{i}\partial_{z} \theta}{(\partial_{z} \theta)^{2}} + \frac{\partial_{i}
\theta \partial_{i} \theta\partial_{z}\partial_{z} \theta}{(\partial_{z} \theta)^{3}}\right). $$
On the other hand the flux $\rho_{t}$ is given by pulling back the differential
form $*_{\epsilon} d\theta$ on the product by the map
$x\mapsto (x, h_{t}(x))$ and this gives 
$$  \rho_{t}= \epsilon \partial_{z} \theta +\frac{1}{\partial_{z}\theta} \sum_{i} (\partial_{i} \theta)^{2}.  $$
So $$ \frac{\partial \rho}{\partial t}= \frac{1}{\partial_{z} \theta} \partial_{z}\left(
\epsilon \partial_{z} \theta + \frac{\sum_{i} (\partial_{i} \theta)^{2}}{\partial_{z} \theta}\right). $$
This is
$$ \partial_{t} \rho_{t}= \epsilon \frac{\partial_{z}\partial_{z} \theta}{\partial_{z} \theta} + 2 \frac{\sum_{i} \partial_{i} \theta \partial_{i} \partial_{z} \theta}{(\partial_{z}
\theta)^{2}} - \frac{\sum_{i}( \partial_{i}\theta)^{2} \partial_{z} \partial_{z} \theta}{(\partial_{z}
\theta)^{3}}. $$
So we see that $\partial_{t} \rho_{t} = -\Delta_{X} h_{t}$, since
 $\epsilon \partial_{z}\partial_{z} \theta = - \sum_{i} \partial_{i}
\partial_{i} \theta$.

Now the normalisation (13) implies that there is a function $\phi_{0}$ on $X$
such that $\rho_{0}= 1-\Delta_{X} \phi_{0}$. For $t>0$  we define $\phi_{t}$ by

 $$\phi_{t}=\phi_{0} +\int_{0}^{t} h_{\tau} d\tau. $$
 We can also regard this family of functions as a single function $\Phi$
 on $X\times[0,1]$,
 Then (15) implies that 
 $\rho_{t}= 1-\Delta_{X}\phi_{t}$ for each $t$. 
    We have
   $$ \frac{\partial^{2} \Phi}{\partial t^{2}}=  \partial_{t} h = \frac{1}{\partial_{z} \theta} $$ and
$$  1-\Delta_{X} \Phi= \epsilon \partial_{z} \theta + \frac{1}{\partial_{z}
\theta} \sum_{i} (\partial_{i} \theta)^{2}. $$
So $$ \frac{\partial^{2} \Phi}{\partial t^{2}}(1-\Delta_{X}\Phi) = \epsilon
+ \sum_{i} \left( \frac{\partial_{i} \theta}{\partial_{z} \theta}\right)^{2}.$$
Now since $$\partial_{i} \partial_{t}\Phi= \partial_{i} h_{t}= \frac{-1}{\partial_{z}
\theta} \partial_{i} \theta$$ we can write the above as
$$ \frac{\partial^{2}\Phi}{\partial t^{2}}(1-\Delta_{X} \Phi) = \epsilon
+ \vert \nabla_{X}\frac{\partial}{\partial t}\Phi \vert^{2} $$ as required.

\item  $\Phi$-equation\ $\Longrightarrow$\ $U$-equation.

\

Here we suppose we have a solution $\Phi(x,t)$ of the $\Phi$ equation and we
write
$\Phi(x,0)=\phi_{0}, \Phi(x,1)=\phi_{1}$.
We essentially take the Legendre transform in the $t$-variable. The discussion is slightly more complicated when $\epsilon=0$, so for simplicity we treat the case when $\epsilon>0$
and $\partial_{t}^{2} \Phi$ is strictly positive. Write
$H_{1}(x), H_{2}(x)$ for the derivatives $\partial_{t} \Phi$ evaluated at
$(x,0),(x,1)$ respectively,  so $H_{0}<H_{1}$. We calculate first in the
open set $\Omega_{H_{0}, H_{1}}$. For each fixed $x\in
X$
and each $z$ in the interval $(H_{0}(x), H_{1}(x))$
there is a $t=t(x,z)$ such that
$ z= \partial_{t} \Phi$. We set
$$  U(x,z)= \Phi(x,0)-\Phi(x,t)+ z t. $$
This defines a function $U$ in $\Omega_{H_{0}, H_{1}}$. We define $U$ outside
this set by setting $U(x,z)=0$ if $z\leq H_{0}(x)$ and $U(x,z)=L(x,z)= \phi_{0}-\phi_{1}-
z$ if $z\geq H_{1}(x)$. It follows from the definitions that $U$ is $C^{1}$,
that $U\geq L$ and that the set where $U>L$ is exactly $\Omega_{H_{0}, H_{1}}$.
We calculate on this set.  
Then
$ \partial_{z} U = t$ and \begin{equation}
\partial_{z}^{2} U = ( \partial_{t}^{2} \Phi)^{-1}.\label{eq:hessinv}\end{equation}
Differentiating with respect to the parameters $x_{i}$ we have
$$  \partial_{i} U= \partial_{i} \phi_{0} - \partial_{i} \Phi, $$
and
$$  \partial_{i}^{2} U= \partial_{i}^{2} \phi_{0}- \partial_{i}^{2} \Phi
- \frac{\partial t}{\partial x_{i}}\frac{ \partial^{2}\Phi}{\partial t\partial
x_{i}}. $$
Differentiating the identity $z=\partial_{t} \Phi$ gives
$$  0= \frac{\partial^{2} \Phi}{\partial t \partial x_{i}} + \frac{\partial^{2}
\Phi}{\partial t^{2}} \frac{\partial t}{\partial x_{i}}, $$
so we can write
$$  \partial_{i}^{2} U= \partial_{i}^{2} \phi_{0}- \partial_{i}^{2} \Phi
+ \frac{1}{\partial_{t}^{2} \Phi} ( \partial_{t} \partial_{i} \Phi)^{2}.
$$Summing over $i$ and using the formula (18) for $\partial_{z}^{2} U$ we obtain
$$ \epsilon \partial_{z}^{2}U - \Delta_{X}U=\frac{1}{\partial_{t}^{2}
\Phi}\left( \epsilon+ \vert \nabla_{Xt} \Phi\vert^{2} \right)-\Delta_{X}\phi_{0}+1,
$$
and so
$$ \Delta_{\epsilon} U=
1-\Delta_{X} \phi_{0}. $$

\

\item $U$-equation $\Longrightarrow$\ $\theta$-equation

Now suppose we have a solution $U$ of $\Delta_{\epsilon}U=\rho_{0}$ in a
domain $\Omega_{H_{0}, H_{1}}$, satisfying the appropriate boundary conditions,
where $\rho_{0}=1-\Delta_{X} \phi_{0}$.
We set 
$$ \theta= \frac{\partial U}{\partial z}. $$
Then $\Delta_{\epsilon} \theta=0$ and $\theta=0,1$ on the two boundary components.
We have to check that the fluxes of $*_{\epsilon}d\theta$ on the boundary components
are $\rho_{i}=1-\Delta_{X}\phi_{i}$. Consider first the boundary component where $z=H_{0}$.
The flux is 
$$  \epsilon \partial_{z} \theta + \frac{\vert \nabla_{X} \theta\vert^{2} }{\partial_{z} \theta}= \epsilon \partial_{z}^{2} F + \frac{1}{\partial_{z}^{2}
F}\sum (\partial_{z}
\partial_{i} F )^{2}. $$
Now we have identities
$$      (\partial_{i} F)(x,H_{0}(x))= 0\ , \ (\partial_{z}F)(x,H_{0}(x))=0.
$$
Differentiating the first of these with repect to $x_{i}$ we get
$$   \partial_{i}^{2}F + \partial_{i} H_{0}\partial_{i} \partial_{z} F=0,
$$ on the boundary. Differentiating the second gives
$$ \partial_{i} \partial_{z} F + \partial_{i} H_{0}\partial_{z}^{2}F =0$$
on the boundary. Combining these we have
$$   (\partial_{z} \partial_{i} F)^{2}= (\partial_{z}^{2} F)(\partial_{i}^{2}
F).$$ Hence the flux is $$\epsilon \partial_{z}^{2} F + \sum_{i} \partial_{i}^{2}
F = \rho_{0}.$$
The argument for the other boundary component $\{z=H_{1}(x)\}$ is similar. 
 \end{itemize}
 \section{Existence results and discussion}
 
 We have set up a class of PDE problems associated to any compact Riemannian
 manifold, and seen that these have three equivalent formulations. In this
 section we will make some remarks about existence results, and comparison
 with the free-boundary literature. This discussion is unfortunately rather
 incomplete, mainly due to the authors limited grasp of the background.
 
 \subsection{Monge-Amp\`ere and the results of Chen}
 
 For a function $\Phi$ on $X\times
 (0,1)$ write $q(\Phi)$ for the nonlinear differential operator
 $$  q(\Phi) = 
  \partial_{t}^{2} \Phi (1-\Delta_{X} \Phi) - \vert \nabla_{X}\frac{\partial}{\partial
  t} \Phi \vert^{2}.
  $$
  So our \lq\lq $\Phi$-equation'' is $q(\Phi)=\epsilon$. 
 When $X$ has dimension $1$---a circle with local coordinate $x$--- we can write $\Delta_{X} = - \partial_{x}^{2}$
 and the equation is the real Monge-Amp\`ere operator
 $$ q(\Phi)_=  \det\left(\begin{array}{cc}  \partial_{t}^{2}\Phi &\partial_{x}\partial_{t}\Phi\\ \partial_{x}\partial_{t}\Phi &1+\partial_{x}^{2} \Phi\end{array} \right) $$

When $X$ has dimension $2$ the operator can be expressed
as a {\it complex} Monge-Amp\`ere operator. That is, we regard $X$ as a Riemann surface
and identify the Laplace operator on $X$ with $i\overline{\partial}\partial$. We
take the product with a circle, with co-ordinate $\alpha$, and let $\tau=t+i\alpha$
be a complex coordinate on the Riemannn surface $S^{1} \times (0,1)$. Then,
in differential form notation, our nonlinear operator  is given by
$$   (\omega_{0}+i\overline{\partial}\partial \Phi)^{2} =q(\Phi) \omega_{0} d\tau
d\overline{\tau}, $$
where $\omega_{0}$ is the Riemannian area form of $X$ lifted to $X\times
S^{1} \times(0,1)$.  Our
Dirichlet problem becomes a Dirichlet problem for $S^{1}$-invariant solutions
of this complex Monge-Amp\`ere equation on $X\times S^{1} \times(0,1)$. 
This was studied by Chen \cite{kn:Chen} and it follows from his results that,
for any $\epsilon>0$ there is a unique solution to our problem, and hence
an affirmative answer to Question 1 in this case. (Chen does not state this
result explicitly, but it follows from the continuity method developed in \cite{kn:Chen},
Section 3, that for any strictly positive smooth function $\nu$ on
$X\times [0,1]$ there is a solution of the equation $q(\Phi)=\nu$ with prescribed
boundary values $\phi_{0}, \phi_{1}$.)  
 
It seems quite likely that the techniques used by Chen can be extended to
the higher dimensional case. The foundation for this should be provided by
a convexity property of the nonlinear operator which we will now derive.
Let $A$ be a symmetric $(n+1)\times (n+1)$ matrix with entries $A_{ij}$
$0\leq i,j\leq n$. Define
$$  Q(A)= A_{00} \sum_{i=1}^{n} A_{ii} - \sum_{i=1}^{n} A_{i0}^{2}. $$
Thus $Q$ is a quadratic function on the vector space of symmetric $(n+1)\times
(n+1)$ matrices.
\begin{lem}
\begin{enumerate}
\item If $A>0$  then $Q(A)>0$ and if $A\geq 0$ then $Q(A)\geq 0$.
\item If $A, B$ are matrices with $Q(A)=Q(B)>0$ and if the entries
$A_{00}, B_{00}$ are positive then for each $s\in [0,1]$
$$  Q( sA+(1-s)B)\geq Q(A)\ ,\  Q(A-B)<0. $$
Moreover, if $A\neq B$ then strict inequality holds. 
\end{enumerate}
\end{lem}
To see the first item, observe that we can change basis in $\bR^{n}\subset
\bR^{n+1}$ to reduce to the case when the block $A_{ij}, 1\leq i,j\leq n$
is diagonal, with entries $b_{i}$ say. Then if $A\geq 0$ we have $A_{00}b_{i}\geq
A_{0i}^{2}$ and so
$$  Q(A)= \sum A_{00}\sum b_{i} - \sum A_{0i}^{2}\geq 0, $$
with strict inequlality if $A>0$. 

For the second item, we just have to observe that $Q$ is induced from a a quadratic form
of  Lorentzian signature on $\bR^{n+2}$ by the linear map
$$\pi: A\mapsto (A_{00}, \sum_{i=1}^{n} A_{ii}, A_{0i}). $$ 
The hypotheses imply that $\pi(A)$ and $\pi(B)$ are in the same component
of a hyperboloid defined by this Lorentzian form and the statements follow
immediately from elementary geometry of Lorentz space.

Using this Lemma we can deduce the uniqueness of the solution to our Dirichlet
problem, in any dimension.
\begin{prop}
If $\phi_{0}, \phi_{1}\in {\cal H}$ then there is at most one solution
$\Phi$ of the equation $Q(\Phi)=\epsilon$ on $X\times [0,1]$ with $1-\Delta_{X} \Phi>0$ for all $t$ and with
$\Phi(x,0)=\phi_{0}(x), \Phi(x,1)=\phi_{1}(x)$.
\end{prop}
We show that the functional $E(\Phi)$ given by (10) is convex with respect to the obvious
linear structure. Thus we consider a $1$-parameter family $\Phi_{s}=\Phi+
s \psi$, with the fixed end points. We have
$$  \frac{d}{ds} E(\Phi_{s})= \int_{0}^{1} \int_{X} 2 \dot{ \Phi_{s}} 
 \dot{\psi} (1-\Delta_{X} \Phi_{s}) -\dot{\Phi}^{2} \Delta_{X} \psi. $$
 Integrating by parts (just as in the derivation of the geodesic equation)
 we obtain

$$  \frac{d}{ds} E(\Phi_{s})= \int_{0}^{1}\int_{X}( q(\Phi_{s})-\epsilon)
\ \psi\ d\mu.$$
 Suppose that $\Phi_{0}, \Phi_{1}$ are two different solutions, so when $s=0,1$ the
 term $q(\Phi_{s})-\epsilon$ in the above expression vanishes pointwise. Item (2) in the lemma
 above implies that for $s\in(0,1)$ we have $q(\Phi_{s})-\epsilon\geq 0$, with strict
 inequality somewhere. This means that $E(\Phi_{1})>E(\Phi_{0})$. Interchanging
 the roles of $\Phi_{0}, \Phi_{1}$ we obtain the reverse inequality, and
 hence a contradiction.

 One can also prove this uniqueness using the maximum principle. Note too
 that the uniqueness is what one would expect, formally, from the negative
 curvature of the space ${\cal H}$ and the convexity of the functional $V$.
 
 \subsection{Comparison with the free-boundary literature}
 
 The author is not at all competent to make this comparison properly. Suffice it to
 say, first, that the problem we are considering is very close to those which
 have been studied extensively in the applied literature. For example, in
 the $\theta$-formulation,  the
 condition of prescribing the pull-back of the flux on the free boundary
 is the same as that in the classical problem of the \lq\lq porous dam''
(\cite{kn:BC} Chapter 8, \cite{kn:EO} Chapter 4.4),  but with the difference that in that case $\rho$
is constant and there are additional boundary conditions on other boundary
components. Second, 
the constructions we have introduced in Section 3 above all appear in this
literature. The transformation from $\theta$ to $\Phi$ taking the harmonic function $\theta$
as a new independent variable is called in (\cite{kn:Crank}, Chapter 5) the \lq\lq isothermal
migration method''. The transformation from the formulation in terms of $\theta$ to that in
terms of $U$ is known as the Baiocchi transformation \cite{kn:BC}, \cite{kn:EO},
\cite{kn:Crank}. The transformation
of the free boundary problem for a linear equation to a nonlinear Dirichlet problem is used in \cite{kn:KN} to derive fundamental regularity results.
 
 An important feature of the $U$-formulation is that it admits a variational
 description. Recall that we are given a function $L= \max(\phi_{0}-\phi_{1}+z,0)$
 on $X\times \bR$ and we seek a $C^{1}$ function $U$ with $U\geq L$ satisfying
 the equation $ \Delta_{\epsilon} U= \rho_{0}$ on the set where $U>L$. This
 can be formulated as follows. We fix a large positive $M$ and consider the
 functional 
 $$   {\cal E}_{M} (U)= \int \frac{1}{2} \vert \nabla_{X}U\vert^{2}+ \epsilon \vert \partial_{z} U \vert ^{2}
 -\rho_{0} U \ d\mu dz , $$
 over the space of functions satisfying the constraint $U\geq L$, where the
 integral is taken over  $X\times [-M,M]$ in $X\times
 \bR$ (which, {\it a posteriori}, should contain the set $\Omega_{H_{0},H_{1}}$
 on which $U>L$). Then the solution minimises ${\cal E}_{M}$ over all functions
 $U\geq L$. This can be used to give another proof of the uniqueness of the
 solution to our problem. It seems likely that it could also be made the basis of an existence
 proof, following standard techniques in the free boundary literature.
 Now recall that our $\Phi$-formulation was based on a variational principle, with Lagrangian (10). To relate the two, we consider any function $\Phi$  on
$X\times [0,1]$ with $\partial_{t} \partial_{t} \Phi\geq 0$ and define $U$
by
the recipe of Section 3. We suppose that $-M< \partial_{t}\Phi(x,0)$ and
$ \partial_{t}\Phi(x,1)<M$
for all $x\in X$. Then we have
\begin{prop}
The functional ${\cal E}_{M}(U)$ is
$$ E(\Phi) + M\int_{X} (1-\Delta_{X}\phi_{0})
(\phi_{0}-\phi_{1})+\frac{1}{2}  \vert \nabla (\phi_{1}-\phi_{0})\vert^{2}
d\mu + (\frac{M^{2}}{2}+\epsilon M)\int_{X} d\mu -\epsilon\int_{X} \phi_{1} d\mu $$
\end{prop}
Thus if we fix $M$ and the end points $ \phi_{0}, \phi_{1}$ the two functionals differ by a constant. The
central step in the proof is the fact that the integrals
$$  \int_{0}^{1}\int_{X} \partial_{t}^{2} \Phi \vert \nabla_{X}\Phi \vert^{2}
\ d\mu\ dt$$
$$\int_{0}^{1} \int_{X} \Delta_{X} \Phi (\partial_{t} \Phi)^{2}\ d\mu\ dt$$ are equal
modulo boundary terms.We leave the full calculation as an exercise for the
reader.

\subsection{The degenerate case}
 So far, in this section, we have discussed the case when $\epsilon>0$. In
 that case the equations we are studying are elliptic. The degenerate case, when $\epsilon=0$,
 is much more delicate. In fact Chen's main concern in \cite{kn:Chen} was to obtain results
 about this case, taking the limit as $\epsilon$ tends to $0$. Chen shows
 that   the Dirichlet problem for $\Phi$, with $\epsilon=0$, has a $C^{1,1}$
 solution
 but the question of smoothness is open. The formulation
 of the problem in terms of the function $U$ has particular advantages here,
 because the problem is set-up as a family of elliptic problems, and the
 issue becomes one of smooth dependence on parameters. (This is related to
 another approach, involving families of holomorphic maps, discussed in \cite{kn:S},
 \cite{kn:D4}.)
 We can express the central question as follows.
  Suppose we have a smooth function
 $\lambda$ on a compact Riemannian manifold $X$ and fix a smooth positive function $\rho$. Let $J$ be the functional
 $$  J(u)= \int_{X}\frac{1}{2} \vert \nabla u \vert^{2} -\rho u. $$
 For each $z\in \bR$ we set $\lambda_{z}= \max(\lambda,z)$ and minimise the functional
 $J$ over the set of functions $u\geq \lambda_{z}$. Suppose we know that there
 is a minimiser $u_{z}$ which is smooth on the open set $\Omega_{z} \subset
 X$ where $u_{z}>\lambda_{z}$. Let $\Omega=\{ (x,z): x\in \Omega_{z}\} \subset X \times \bR$.
\begin{Question}
In this situation, does $u_{z}$ vary smoothly with $z$ in $\Omega$?
\end{Question}   
The interesting case here seems to be when $z$ is a critical value of $g$.  
  \section{Relation with Nahm's equations}
  
  We recall that Nahm's equations are a system of ODE for three functions
  $T_{1},T_{2}, T_{3}$ taking values in a fixed Lie algebra:
  \begin{equation} \frac{dT_{i}}{dt}= [T_{j},T_{k}], \label{eq:Nahm1} \end{equation}
  where $i,j,k$ run over cyclic permutations of $1,2,3$. To simplify notation,
  let us fix on the Lie algebra $u(n)$.  It is equivalent
  (at least in the finite-dimensional case) to introduce a fourth function
  $T_{0}$ and consider the equations
  \begin{equation}  \frac{dT_{i}}{dt}+ [T_{0}, T_{i}]= [T_{j},T_{k}], \label{eq:Nahm2}
  \end{equation}
  with the action of the \lq\lq gauge group'' of $U(n)$-valued functions $u(t)$:
  $$   T_{i}\mapsto u T_{i} u^{-1}, T_{0}\mapsto u T_{0} u^{-1}- \frac{du}{dt} u^{-1}  $$
which preserved solutions to (20). (That is, using the gauge group we can transform
$T_{0}$ to $0$.) The equations imply that
 \begin{equation}  \frac{d}{dt} (T_{2}+iT_{3}) = [T_{0}+i T_{1}, T_{2}+ iT_{3}],\label{eq:invar}\end{equation} so 
   $T_{2}+iT_{3}$ moves in a single adjoint orbit in the Lie algebra of $GL(n,\bC)$.
   Conversely if we fix some $B$ in this complex Lie algebra,  introduce a function $g(t)$ taking values in $GL(n,\bC)$
   and define skew-Hermitian matrices $T_{i}(t)$ by
   $$   T_{0}+iT_{1}= \frac{dg}{ds} g^{-1},$$
   $$  T_{2}+iT_{3}= g B g^{-1}, $$
   then two of the three Nahm equations are satisfied identically. The remaining
   equation can be expressed in terms of the function $h(t)= g^{*}(t)g(t)$,
   taking values in the space ${\cal H}$ of positive definite Hermitian matrices,
   which we can also regard as the quotient space $GL(n,\bC)/U(n)$. This
   equation for $h(t)$ is a second order ODE which is the Euler-Langrange
   equation for the Lagrangian
   $$   E(h)= \int \vert \frac{dh}{dt} \vert^{2}_{{\cal H}} + V_{B}(h)  dt. $$
   Here $\vert\ \vert_{{\cal H}}$ denotes the standard Riemannian metric
   on ${\cal H}$. The function $V$  on ${\cal H}$ is
   $$  V_{B}(h)= {\rm Tr} (hBh^{-1}B^{*}). $$
   If $g$ is any element of $GL(n,\bC)$ with $g^{*}=h$ then
   $$   V_{B}(h)= \vert g B g^{-1} \vert^{2}, $$
   so $V_{B}$ is determined by the norm of matrices in the adjoint orbit
   of $B$. (See \cite{kn:D1} for details of the manipulations involved in
  all the above.) The result in \cite{kn:D1}, mentioned in the introduction to this article,
    is that for any two points $h_{0}, h_{1}\in {\cal H}$ there is a unique
    solution $h(t)$ to the Euler-Lagrange equations for $t\in [0,1]$ with
    $h(0)=h_{0}, h(1)=h_{1}$.

    These constructions go over immediately to the case when $U(n)$ is replaced
    by any compact Lie group and $GL(n,\bC)$ by the complexified group. We
    want to extend them to the situation where $U(n)$ is replaced by the
    group ${\cal G}$ of Hamiltionian diffeomorphisms of a surface $\Sigma$ with a fixed
    area form (or more precisely, the extension of this group given by a
    choice of Hamiltionian). The essential difficulty is that this group
    does not have a complexification. However, as explained in \cite{kn:D3},
    \cite{kn:M}, \cite{kn:S},  the
    space ${\cal H}$ of Kahler potentials behaves formally like the quotient
    space ${\cal G}^{c}/{\cal G}$ for a fictitious group ${\cal G}^{c}$.
    Thus the problem we have formulated in Section 2 can be viewed as an
    analogue of the desired kind provided that our potential function $V$
    can be seen as an analogue of $V_{B}$ in the finite-dimensional case.

    If we have a path $\phi_{t}$ in ${\cal H}$ with $\phi_{0}=0$ and a function $\beta:\Sigma
    \rightarrow \bC$ we can write down a differential equation
    for a one-parameter family $\beta_{t}$ which corresponds, formally, to
    the adjoint action of the complexified group ${\cal G}^{c}$, with the
    initial condition $\beta_{0}=\beta$. The equation has the shape
    $$   \frac{\partial \beta_{t}}{\partial t}= \nabla \dot{\phi} \overline{\partial}
    \beta_{t}. $$
    The problem is that this evolution equation will not have solutions,
    even for a short time, in general. But if we suspend for a moment our
    assumption that we are working over  a compact Riemann surface and suppose that $\beta$ is a {\it
    holomorphic} function then there is a trivial solution $\beta_{t}=\beta$.
    So, formally, the functional $V_{\beta}$ on ${\cal H}$ is given by the
    $L^{2}$ norm of $\beta$ with respect to the meaure $d\mu_{\phi}$:
    $$   \int (1-\Delta_{X} \phi) \vert \beta \vert^{2}. $$
    Even if this integral is divergent, the variation with respect to compactly
    supported variations in $\phi$ is well-defined, and this is what corresponds
    to the gradient of $V_{B}$ appearing in the equations of motion. Moreover,
    we can integrate by parts to get another formal representation of a functional
    with the same variation
    $$  -\int \phi \Delta_{X} \vert \beta \vert^{2}=  \int \phi \vert \nabla
    \beta \vert^{2}. $$
    
    Now take the compact Riemann surface $\Sigma$ to be a 2-torus,
    and identify the space ${\cal H}$ with periodic Kahler potentials on
    the universal cover $\bC$. On this cover the identity function $\beta$
    is holomorphic, and we see from the above that the formal expression
    $$  V_{\beta}= \int_{\bC} \phi, $$
    is analogous to the function $V_{B}$ in the finite-dimensional case.
    Of course the integrand is periodic and so the integral will be divergent
    but we can return to the compact surface $\Sigma$ and consider the
    well-defined functional
    $$  V_{\beta}(\phi) = \int_{\Sigma} \phi $$
    which will generate the same equations of motion. So we see that, modulo
    some blurring of the distinction between $\Sigma$ and its universal cover,
    the functional we have been considering is indeed analogous to that in
    the finite-dimensional case.
    
    Using the transformation from the $\Phi$ equation to the $\theta$ equation,
    we obtain a 
     relation between  Nahm's equations for the Hamiltonian diffeomorphisms
    of a surface  and 
     harmonic functions on $\bR^{3}$. This  can be seen in other ways. Most
     directly, we consider three one-parameter families of functions $h_{i}(t)$ on a surface
     $\Sigma$ with an area form which satisfy:
     \begin{equation}  \frac{dh_{i}}{dt}= \{ h_{j}, h_{k}\}, \label{eq:Poisson}\end{equation}
     where $\{\ ,\ \}$ is the Poisson bracket. We think of these as a one-parameter
     family of maps $\underline{h}_{t}:\Sigma \rightarrow \bR^{3}$, and assume
     for simplicity that these are embeddings, with disjoint images. Then
     it is a simple exercise to show that the equations (22) imply that the
     images $\underline{h}_{t}(\Sigma)$ are the level sets of a harmonic
     function on a domain in $\bR^{3}$. From another point of view, the geometric
     structure defined by a solution to the $\Phi$ equation is an $S^{1}$
     invariant Kahler metric $\Omega=\omega_{0}+ i \overline{\partial} \partial \Phi$  on $\Sigma\times S^{1} \times (0,1)$ with volume form
$$ \Omega^{2} = d\tau d\beta \overline{d\tau d\beta}. $$
Since $d\tau d\beta$ is an $S^{1}$-invariant holomorphic $2$-form, what we
have is an $S^{1}$ invariant {\it hyperkahler} structure. Then the relation
with harmonic functions appears as  the Gibbons-Hawking construction for hyperkahler
metrics.

The development above is rather limited, since we have only been able formulate
an analogue of our Nahm's equation problem for a single function $\beta$.
One
can go further, and arrive at other interesting free boundary problems. Consider
for example the case when the surface $\Sigma$  is the $2$-sphere with the
standard area form, and the orientation-reversing map $\sigma:\Sigma\rightarrow \Sigma$ given by reflection in the $x_{1}, x_{2}$ plane. Now consider maps
$\beta: \Sigma \rightarrow \bC$ with $\beta= \beta \circ \sigma$ which are
diffeomorphisms on each hemisphere. Then the push-forward of the area form
on the upper hemisphere defines a  $2$-form $\rho_{\beta}$ on $\bC$ with support
in a topological disc $\beta(\Sigma)\subset \bC$. (The form $\rho_{\beta}$ will not usually be
smooth, but will behave like $d^{-1/2}$ where $d$ is the distance to the
boundary of $\beta(\Sigma)$.) Clearly the form $\rho_{\beta}$ determines
$\beta$ up to the action of the $\sigma$-equivariant Hamiltionian diffeomorphisms
of $\Sigma$. Suppose that $h$ is a $\sigma$-invariant function on $\Sigma$.
We can regard this as an element of the Lie algebra of ${\cal G}^{c}$ and
consider its action on $\beta$. This is given by $\Delta_{\bC} h$ where
$h$ is thought of as a function on $\bC$, vanishing outside $\beta(\Sigma)$.
So a reasonable candidate for a model of the quotient of the space of maps $\beta$ by the action
of ${\cal G}^{c}$ is given by the following. We consider $2$-forms $\rho$
supported on topological discs in $\bC$, with singularities at the boundary
of the kind arising above, and impose the equivalence relation that
$\rho_{0}\sim \rho_{1}$ if there is a compactly supported harmonic function
$F$ on $\bC$ with $\Delta F = \rho_{0}- \rho_{1}$.

Now let $\theta(x_{1}, x_{2}, z)$ be a harmonic function on an open set  $\Omega\subset \bR^{3}$, with
$\theta(x_{1}, x_{2}, z)= \theta(x_{1}, x_{2}, -z)$. Suppose that $\Omega$
is diffeomorphic to $S^{2} \times (0,1)$, that $\theta=0$  on the
inner boundary component $\Sigma_{0}$ and $\theta=1$ on the outer
boundary component $\Sigma_{1} $. Suppose also that the projections
of $\Sigma_{0}, \Sigma_{1} $ to the $(x_{1}, x_{2})$ plane
are diffeomorphisms on each upper hemisphere, mapping to a pair of topological
disc $D_{0}\subset D_{1}$.  Then the flux of $\nabla \theta$
on each boundary component pushes forward to define a pair of compactly supported
$2$-forms $\rho_{0}, \rho_{1}$ on $\bC$. These are equivalent in the sense
above, since $\rho_{0}-\rho_{1}=\Delta_{\bC} F$ for the function
$$   F(x_{1}, x_{2}) =\int z \frac{\partial \theta}{\partial z} dz, $$
where the integral is taken over the intersection of the vertical line through
$(x_{1}, x_{2}, 0)$ with $\Omega$. Our hypotheses imply that $F\geq 0$, and
$F$ is supported on the larger disc $D_{1}$. 

The question we are lead to is the following
\begin{Question}
Suppose $D_{0}\subset D_{1}$ are topological discs in $\bC$, that $\rho_{i}$
are
$2$-forms supported on $D_{i}$ and that there is a non-negative function $F$ on $\bC$, supported
on $D_{1}$,  with $\rho_{0}-\rho_{1} =\Delta_{\bC} F$ (where the Laplacian
is defined in the distributional sense). Do $\rho_{0}, \rho_{1}$ arise from
a  unique harmonic function
$\theta$ on a domain in $\bR^{3}$, by the construction above? 
\end{Question}
(For simplicity we have not specified precisely what singularities should
be allowed in the forms $\rho_{i}$: this specification should be a part of
the question.)
  
   Hitchin showed in \cite{kn:H} that Nahm's equations form an integrable
   system. The root of this is the invariance of the conjugacy class given
   by (21), together with the family of similar statements that arise from    the $SO(3)$ action on the set-up. In this vein, we can write down infinitely
   many conserved quantities for the solutions of our  equation (11) on
   the Riemannian manifold ${\cal H}$. Let $f_{\lambda}$ be an eigenfunction
   of the Laplacian $\Delta_{X}$, with eigenvalue $\lambda>0$. Then we have
   \begin{prop}
   
   For any $\epsilon>0$, if $\phi_{t}$ satisfies (11) then the quantity
   $$     \int_{X} \exp\left(\sqrt{\frac{\lambda}{\epsilon}} \ \dot{\phi}\right)f_{\lambda}
   (1-\Delta \phi) \ d\mu, $$
   does not vary with $t$.
   \end{prop}
   
   This becomes rather transparent in the $\theta$-formulation, using the
   fact that the function 
   $$   K_{\lambda}(x,z)= f_{\lambda}(x) \exp\left( \sqrt{\frac{\lambda}{\epsilon}}
   \ z \right), $$
   satisfies $\Delta_{\epsilon} K_{\lambda}=0$.

     \
     
     The author is grateful to Professors Colin Atkinson, Xiuxiong Chen, Darryl Holm
     and John Ockendon for helpful discussions.



\end{document}